# A Brief Survey of the History of the Calculus of Variations and its Applications

James Ferguson

jcf@uvic.ca

University of Victoria

## Abstract

In this paper, we trace the development of the theory of the calculus of variations. From its roots in the work of Greek thinkers and continuing through to the Renaissance, we see that advances in physics serve as a catalyst for developments in the mathematical theory. From the 18th century onwards, the task of establishing a rigourous framework of the calculus of variations is studied, culminating in Hilbert's work on the Dirichlet problem and the development of optimal control theory. Finally, we make a brief tour of some applications of the theory to diverse problems.

**Introduction**

Consider the following three problems:

1) What plane curve connecting two given points has the shortest length?

2) Given two points A and B in a vertical plane, find the path AMB which the movable particle M will traverse in shortest time, assuming that its acceleration is due only to gravity.

3) Find the minimum surface of revolution passing through two given fixed points, $(x_A, y_A)$ and $(x_B, y_B)$.

All three of these problems can be solved by the calculus of variations. A field developed primarily in the eighteenth and nineteenth centuries, the calculus of variations has been applied to a myriad of physical and mathematical problems since its inception. In a sense, it is a generalization of calculus. Essentially, the goal is to find a path, curve, or surface for which a given function has a stationary value. In our three introductory problems, for instance, this stationary value corresponds to a minimum.

The variety and diversity of the theory's practical applications is quite astonishing. From soap bubbles to the construction of an ideal column and from quantum field theory to softer spacecraft landings, this venerable branch of mathematics has a rich history and continues to spring upon us new surprises. Its development has also served as a catalyst for theoretical advances in seemingly disparate fields of mathematics, such as analysis, topology, and partial differential equations. In fact, at least two modern (i.e. since the beginning of the twentieth century) areas of research can claim the calculus of variations as a common ancestor; namely Morse theory and optimal control theory. Since the theory was initially developed to tackle physical problems, it is not surprising that variational methods are at the heart of modern approaches to problems in theoretical physics. More surprising is that the calculus of variations has been applied to problems in economics, literature, and interior design!

In the course of this paper, we will trace the historical development of the calculus of variations. Along the way, we will explore a few of the more interesting historical problems and applications, and we shall highlight some of the major contributors to the theory. First, let us get an intuitive sense of the theory of the calculus of variations with the following mathematical interlude, which might be found along similar lines in an applied math or physics text (e.g. [2] and [5]).

**Mathematical Background**

In this section we derive the differential equation that $y(x)$ must obey in order to minimize the integral

$$I = \int_{x_A}^{x_B} f(x, y, y')dx$$

where $x_A$, $x_B$, $y(x_A) = y_A$, $y(x_B) = y_B$ and $f$ are all given, and $f$ is assumed to be a twice-differentiable function of all its arguments. Let us denote the function which minimizes $I$ to be $y(x)$. Now consider the one-parameter family of comparison functions (or test functions), $\tilde{y}(x, \varepsilon)$, which satisfy the conditions:

a) $\tilde{y}(x_A, \varepsilon) = y_A$, $\tilde{y}(x_B, \varepsilon) = y_B$ for all $\varepsilon$;
b) $\tilde{y}(x, 0) = y(x)$, the desired minimizing function;
c) $\tilde{y}(x, \varepsilon)$ and all its derivatives through second order are continuous functions of $x$ and $\varepsilon$.

For a given comparison function, the integral

$$I(\varepsilon) = \int_{x_A}^{x_B} f(x, \tilde{y}, \tilde{y}')dx$$

is clearly a function of $\varepsilon$. Also, since setting $\varepsilon = 0$ corresponds, by condition (b), to replacing $\tilde{y}$ by $y(x)$ and $\tilde{y}'$ by $y'(x)$, we see that $I(\varepsilon)$ should be a minimum with respect to $\varepsilon$ for the value $\varepsilon = 0$ according to the designation that y(x) is the actual minimizing function. This is true for any $\tilde{y}(x, \varepsilon)$.

A necessary condition for a minimum is the vanishing of the first derivative. Thus we have

$$\left[\frac{dI}{d\varepsilon}\right]_{\varepsilon=0} = 0$$

as a necessary condition for the integral to take on a minimum value at $\varepsilon = 0$. Differentiating with respect to $\varepsilon$ (remembering that x is a function only of $y$ and $\tilde{y}$), we get:

$$\frac{dI}{d\varepsilon} = \int_{x_A}^{x_B}\left[\frac{\partial f}{\partial \tilde{y}}\frac{\partial \tilde{y}}{\partial \varepsilon} + \frac{\partial f}{\partial \tilde{y}'}\frac{\partial \tilde{y}'}{\partial \varepsilon}\right]dx$$

and by condition (c), we can write this as:

$$\frac{dI}{d\varepsilon} = \int_{x_A}^{x_B}\left[\frac{\partial f}{\partial \tilde{y}}\frac{d\tilde{y}}{d\varepsilon} + \frac{\partial f}{\partial \tilde{y}'}\frac{d}{dx}\left(\frac{d\tilde{y}'}{d\varepsilon}\right)\right]dx.$$

Integrating the second term by parts gives us:

$$\frac{dI}{d\varepsilon} = \int_{x_A}^{x_B}\frac{\partial f}{\partial \tilde{y}}\frac{d\tilde{y}}{d\varepsilon}dx + \left[\frac{d\tilde{y}}{d\varepsilon}\frac{\partial f}{\partial \tilde{y}'}\right]_{x_A}^{x_B} - \int_{x_A}^{x_B}\frac{d\tilde{y}}{d\varepsilon}\frac{d}{dx}\left(\frac{df}{d\tilde{y}'}\right)dx.$$

Now by condition (a), $\tilde{y}(x_A,\varepsilon) = y_A$ and $\tilde{y}(x_B,\varepsilon) = y_B$ for all $\varepsilon$. Therefore,

$$\left.\frac{d\tilde{y}}{d\varepsilon}\right|_{x=x_A} = 0 = \left.\frac{d\tilde{y}}{d\varepsilon}\right|_{x=x_B}$$

and in the end, we get:

$$\frac{dI}{d\varepsilon} = \int_{x_A}^{x_B}\left[\frac{\partial f}{\partial \tilde{y}} - \frac{d}{dx}\left(\frac{\partial f}{\partial \tilde{y}'}\right)\right]\frac{d\tilde{y}}{d\varepsilon}dx.$$

We now require that $I(\varepsilon)$ have a minimum at $\varepsilon = 0$, that is

$$\left[\frac{dI}{d\varepsilon}\right]_{\varepsilon=0} = \int_{x_A}^{x_B}\left[\frac{\partial f}{\partial \tilde{y}} - \frac{d}{dx}\left(\frac{\partial f}{\partial \tilde{y}'}\right)\right]_{\varepsilon=0}\left[\frac{d\tilde{y}}{d\varepsilon}\right]_{\varepsilon=0}dx.$$

If we set $\varepsilon = 0$, this is the same as setting $\tilde{y}(x,\varepsilon) = y(x)$, $\tilde{y}'(x,\varepsilon) = y'(x)$, and $\tilde{y}''(x,\varepsilon) = y''(x)$. (Note that the integrand depends on $\tilde{y}''$, and in taking the limit $\varepsilon = 0$, we need to know that the second derivative $\tilde{y}''(x,\varepsilon)$ is a continuous function of its two variables. This is guaranteed by condition (c).)

Now if we set

$$\left[\frac{d\tilde{y}}{d\varepsilon}\right]_{\varepsilon=0} = \eta(x),$$

we obtain

$$\int_{x_A}^{x_B}\left[\frac{\partial f}{\partial \tilde{y}} - \frac{d}{dx}\left(\frac{\partial f}{\partial \tilde{y}'}\right)\right]\eta(x)dx = 0.$$

Now $\eta(x)$ vanishes at $x_A$ and $x_B$ by condition (a) and it is continuous and differentiable by condition (c). However, aside from these qualities, $\eta(x)$ is completely arbitrary. Therefore, for the integral above to vanish, we must have

$$\frac{\partial f}{\partial y} - \frac{d}{dx}\left(\frac{\partial f}{\partial y'}\right) = 0.$$

This is known as the Euler-Lagrange equation, which is used to develop the Lagrangian formulation of classical mechanics. If we expand the total derivative with respect to x, we get

$$\frac{\partial f}{\partial y} - \frac{\partial^2 f}{\partial x \partial y} - \frac{\partial^2 f}{\partial y \partial y'} y' - \frac{\partial^2 f}{\partial^2 y'^2} y'' = 0.$$

This is a second-order differential equation, whose solution is a twice-differentiable minimizing function $y(x)$, provided a minimum exists. Note that our initial condition of

$$\left[\frac{dI}{d\varepsilon}\right]_{\varepsilon=0} = 0$$

is only a necessary condition for a minimum. The solution $y(x)$ could also produce a maximum or an inflection point. In other words, $y(x)$ is an extremizing function. [5]

**Hero and the Principle of Least Time**

Probably the first person to seriously consider minimization problems from a scientific point of view was Hero of Alexandria, who lived sometime between 150 BC and 300 AD. He studied the optics of reflection and pointed out, without proof, that reflected light travels in a way that minimizes its travel time. This is a precursor to Fermat's principal of least time.

Hero showed that when a ray of light is reflected by a mirror, the path taken from the object to the observer's eye is shorter than any other possible path so reflected. It is worthwhile to quote from Hero's *Catoptrics*:

> Practically all who have written of dioptrics and of optics have been in doubt as to why rays proceeding from our eyes are reflected by mirrors and why the reflections are at equal angles. Now the proposition that our sight is directed in straight lines proceeding from the organ of vision may be substantiated as follows. For whatever moves with unchanging velocity moves in a straight line... for because of the impelling force the object in motion strives to move over the shortest possible distance, since it has not the time for slower motion, that is, for motion over a longer trajectory. The impelling force does not permit such retardation. And so, by reason of its speed, the object tends to move over the shortest path. But the shortest of all lines having the same end points is a straight line... Now by the same reasoning, that is, by a consideration of the speed of the incidence and the reflection, we shall prove that these rays are reflected at equal angles in the case of plane and spherical mirrors. For our proof must again make use of minimum lines. [20]

**Pappus and Isoperimetric Problems**

Once upon a time, kings would reward exceptional civil servants and military personnel by giving them all the land that they could encompass by a ploughed furrow in a specified period of time. In this way, the problem of finding the plane curve of a given length which encloses the greatest area, or the isoperimetric problem, was born [5]. Pappus of Alexandria (c.290 AD - c.350 A.D.) was not the first person to consider isoperimetric problems. However, in his book Mathematical Collection, he collected and systematized results from many previous mathematicians, drawing upon works from Euclid (325 BC - 265 BC), Archimedes (287 BC - 212 BC), Zenodorus (200 BC - 140 BC), and Hypsicles (190 BC - 120 BC). This topic is often linked to the five so-called Platonic solids (pyramid, cube, octahedron, dodecahedron, and icosahedron).

In Book 5 of the *Mathematical Collection*, Pappus compares figures with equal contours (or surfaces) to see which has the greatest area (or volume). We can summarize the main mathematical contents of Book 5 in the following way:

1. Among plane figures with the same perimeter, the circle has the greatest area.
2. Among solid figures with the same area, the sphere has the greatest volume.
3. There are five and only five regular solids.

Apart from these three primary results, Pappus also notes the following secondary points:

1. Given any two regular plane polygons with the same perimeter, the one with the greater number of angle has the greater area, and consequently,
2. Given a regular plane polygon and a circle with the same perimeter, the circle has the greater area.
3. A circle has the same area as a right-angled triangle whose base is equal to the radius and whose height is equal to the circumference of the circle.
4. Of isoperimetric polygons with the same number of sides, a regular polygon is greater than an irregular one.
5. Given any segments with the same circumference, the semicircle has the greatest area.
6. There are only five regular solid bodies.
7. Given a sphere and any of the five regular solids with equal surface, the sphere is greater.
8. Of solid bodies with the same surface, the one with more faces is the greatest.
9. Every sphere is equal to a cone whose base is the surface of the sphere and whose height is its radius.

Pappus appears to have been a master of demonstrating what had already been shown. In fact, item 9 from the list above was well known to the world as being proved by Archimedes and was even engraved on his tombstone. In spite of this, his works were a useful collection of facts related to problems about isoperimetry [6].

**Fermat (1601 - 1665)**

A more serious and more general minimization problem in optics was studied in the mid-17th century by the French mathematician Pierre de Fermat (1601-1665). He believed that "nature operates by means and ways that are 'easiest and fastest'" but not always on shortest paths. When it came to light rays, Fermat believed that light travelled more slowly in a denser medium. (While this may seem intuitive to us, Descartes believed the opposite - that light travelled *faster* in a denser medium.) He was able to show that the time required for a light ray to traverse a neighbouring virtual path differs from the time actually taken by a quantity of the second order [20].

We can state Fermat's principle mathematically as:
$$\delta \int_P^Q \frac{ds}{v} = 0 \, ,$$
where P and Q are the starting- and end-points of the path, *v* the velocity at any point and *ds* an element of the path. The equation indicates that the variation of the integral is zero, i.e., the difference between this integral taken along the actual path and that taken along a neighbouring path is an infinitesimal quantity of the second order in the distance between the paths.

However, this disagreement with the great René Descartes (1596-1650) was the cause of much personal and public agony for Fermat. One can sense the style and wit of Fermat, in the following excerpt from a letter he wrote to Clerselier (a defender of Descartes) in May 1662:

> I believe that I have often said both to M. de la Chambre and to you that I do not pretend, nor have I ever pretended to be in the inner confidence of Nature. She has obscure and hidden ways which I have never undertaken to penetrate. I would have only offered her a little geometrical aid on the subject of refraction, should she have been in need of it. But since you assure me, Sir, that she can manage her affairs without it, and that she is content to follow the way that has been prescribed to her by M. Descartes, I willingly hand over to you my alleged conquest of physics; and I am satisfied that you allow me to keep my geometrical problem - pure and *in abstracto*, by means of which one can find the path of a thing moving through two different media and seeking to complete its movement as soon as it can [20].

**Newton and Surfaces of Revolution in a Resisting Medium**

The first real problem in the calculus of variations was studied by Sir Isaac Newton (1643-1727) in his famous work on mechanics, *Philosophiae naturalis principia mathematica* (1685), or the *Principia* for short [11]. Newton examined the motion of bodies in a resisting medium. First, he considered a specific case, that of the motion of a frustum of a cone moving through a resisting medium in the direction of its axis. This problem can be solved using the ordinary (i.e. pre-existing) theory of maxima and minima, which Newton showed. Next, Newton considered a more general problem. Suppose a body moves with a given constant velocity through a fluid, and suppose that the body covers a prescribed maximal cross section (orthogonal to the velocity vector) at its rear end. Find the shape of the body which renders its resistance minimal.

This was the first problem in the field to be clearly formulated and also the first to be correctly solved, thus marking the birth of the theory of the calculus of variations. The geometrical technique used by Newton was later adopted by Jacob Bernoulli in his solution of the Brachistochrone and was also later systematized by Euler. Aside from giving birth to an entire field, a further point of interest about Newton's study of motion in a resisting medium is that it is actually one of the most difficult problems ever tackled by variational methods until the twentieth century. Firstly, the formulation of the problem requires several assumptions to be made regarding the resisting medium and the nature of the resistance experienced by the moving body. As it turns out, the restrictions imposed by Newton are only valid for bodies moving at a velocity greater than the speed of sound for the given medium. Secondly, the problem can possess solution curves having a corner (i.e. a discontinuous slope) which, when expressed parametrically, may not have a solution in the ordinary sense [11]. This foreshadows twentieth century developments in optimal control theory.

We should make a few further remarks regarding Newton's solution to this problem, as appeared in the *Principia*. Anyone who has ever been baffled by a mathematical text before will find solace in the fact that Newton's solution appeared without a suggestion or hint as to how to derive it. Furthermore, none of Newton's contemporaries, including Leibniz (but with the possible exception of Huygens), could grasp the fundamental ideas behind Newton's technique. The mathematical community was completely baffled. Eventually, an astronomy professor at Oxford, named David Gregory, persuaded Newton to write out an analysis of the problem for him in 1691. After studying Newton's detailed exposition, Gregory communicated it to his students, and thereby the rest of the world, through his Oxford lectures in the fall of 1691. Since that time, numerous studies have been undertaken involving more general considerations (i.e. more realistic types of resistance, non-symmetric surfaces, etc.). A good overview can be found in [4].

**The Brachistochrone**

The most famous problem in the history of the subject is undoubtedly the problem of the Brachistochrone. In June of 1696, Johann Bernoulli (1667-1748), a member of the most famous mathematical family in history, issued an open challenge to the mathematical world with the following problem (problem (2) from the Introduction above):

> Given two points A and B in a vertical plane, find the path AMB which the movable particle M will traverse in shortest time, assuming that its acceleration is due only to gravity.

The problem is in fact based on a similar problem considered by Galileo Galilei (1564-1642) in 1638. Galileo did not solve the problem explicitly and did not use methods based on the calculus. Due to the incomplete nature of Galileo's work on the subject, Johann was fully justified in bringing the matter to the attention of the world. After stating the problem, Johann assured his readers that the solution to the problem was very useful in mechanics and that it was not a straight line but rather a curve familiar to geometers. He gave the world until the end of 1696 to solve problem, at which time he promised to publish his own solution. At the end of the year, he published the challenge a second time, adding an additional problem (one of a geometrical nature), and extending his deadline until Easter of 1697.

At the time of the initial challenge to the world, Johann Bernoulli had also sent the problem privately to one of the most gifted minds of the day, Gottfried Wilhelm Leibniz (1646-1716), in a letter dated 9 June 1696. A short time later, he received a complete solution in reply, dated 16 June 1696! In our modern society, which has become obsessed doing everything "as soon as possible", focusing so much on speed that we often sacrifice quality, it is refreshing to see that technology is not a prerequisite for timeliness. It also gives us an indication of Leibniz's genius. It was in correspondence between Leibniz and Johann Bernoulli that the name Brachistochrone was born. Leibniz had originally suggested the name *Tachistoptotam* (from the Greek *tachistos*, swiftest, and *piptein*, to fall). However, Bernoulli overruled him and christened the problem under the name *Brachistochrone* (from the Greek *brachistos*, shortest, and *chronos*, time) [11].

The other great mathematical mind of the day, Newton, was also able to solve the problem posed by Johann Bernoulli. As legend has it, on the afternoon of 29 January 1697, Newton found a copy of Johann Bernoulli's challenge waiting for him as he returned home after a long day at work. At this time, Newton was Warden of the London Mint. By four o'clock that morning, after roughly twelve hours of continuous work, Newton had succeeded in solving both of the problems found in Bernoulli's challenge! That same day, Newton communicated his solution anonymously to the Royal Society. While it is quite a feat, comparable to that of Leibniz's rapid response to Bernoulli, one should note that Bernoulli himself claimed that neither problem should take "*a man capable of it more than half an hour's careful thought*." As Ball slyly notes [3], since it actually took Newton twelve hours, it is "*a warning from the past of how administration dulls the mind*." Indeed, it is rather surprising that it took Newton so long, considering the similarities that the Brachistochrone problem has with Newton's previously solved problem of bodies in a resisting medium.

When Johann originally posed the problem, it is likely that his main motivation was to fuel the fire of his bitter feud with elder brother, Jacob Bernoulli (1654-1705). Johann had publicly described his brother Jacob as incompetent and was probably using the Brachistochrone problem, which he has already solved, as a means of publicly triumphing over his brother. Such an attitude towards ones contemporaries prompted one scholar to remark that it must have been Johann Bernoulli who first said the words, "It is not enough for you to succeed; your colleagues must also fail." [5]

In the end, Jacob Bernoulli was able to solve the problem set to him by his brother, joining Leibniz, Newton, and l'Hôpital as the only people to correctly solve the problem. It is interesting to note that even though Newton sent in his result anonymously, Johann Bernoulli was not fooled. He later wrote to a colleague that Newton's unmistakable style was easy to spot and that "he knew the lion from his touch." Far from being gracious, however, Johann was quick to proclaim his superiority over others when summarizing the results of his challenge:

> I have with one blow solved two fundamental problems, one optical and the other mechanical and have accomplished more than I have asked of others: I have shown that the two problems, which arose from totally different fields of mathematics, nevertheless possess the same nature.

Bernoulli refers to the fact that he was the first to publicly demonstrate that the least time principle of Fermat and the least time nature of the Brachistochrone are two manifestations of the same phenomenon.

Let us exhibit a solution for the Brachistochrone problem, not in the geometrical language of the times, but rather in the more modern way that was developed later by Euler and Lagrange (as we shall soon see).

> Let us take *A* as the origin in our coordinate system, assume that the particle of mass *m* has zero initial velocity, and assume that there is no friction. Let us also take the *y*-axis to be directed vertically downward. The speed along the curve AMB is $v = ds/dt$ and thus, the total time of descent is
>
> $$ I = \int_A^B \frac{ds}{v} = \int_{x_A=0}^{x_B} \frac{\sqrt{1+y'^2}}{v} dx. $$
>
> Now we know by conservation of energy that the change in kinetic energy must equal the change in potential energy. Therefore, we can write
>
> $$ \frac{1}{2}mv^2 = mgy, $$
>
> so that the functional to be minimized becomes
>
> $$ I = \frac{1}{\sqrt{2g}} \int_0^{x_B} \sqrt{\frac{1+y'^2}{y}} dx. $$
>
> Now we can use the Euler-Lagrange equation to obtain (neglecting the constant factor of $1/\sqrt{2g}$ )
>
> $$ \frac{y'^2}{\sqrt{y(1+y'^2)}} - \sqrt{\frac{1+y'^2}{y}} = C \quad \text{or} \quad \frac{1}{y(1+y'^2)} = C^2. $$
>
> Setting $1/C^2 = 2a,$ we obtain
>
> $$ y' = \sqrt{\frac{2a-y}{y}}, $$

and integration yields

$$x - x_0 = \int \sqrt{\frac{y}{2a - y}} dy.$$

After making the change of variables $y = a(1 - \cos\theta)$, the integral becomes

$$x - x_0 = 2a \int \sin^2 \frac{\theta}{2} d\theta = a(\theta - \sin\theta).$$

Therefore, the solution to the brachistochrone problem, in parametric form, is

$$x = a(\theta - \sin\theta) + x_0, \qquad y = a(1 - \cos\theta).$$

These are the equations of a cycloid generated by the motion of a fixed point on the circumference of a circle of radius $a$ which rolls on the positive side of the line $y = 0$, that is, on the underside of the $x$-axis. There exists one and only one cycloid through the origin and the point $(x_B, y_B)$; a suitable choice of $a$ and $x_0$ will give this cycloid [5].

## Euler, Maupertuis, and the Principle of Least Action

The brilliant and prolific Swiss mathematician Leonhard Euler (1707-1783) had close ties to the Bernoulli family. Not only was his father, Paul Euler, friends with Johann but Paul had also lived in Jakob's house while he studied theology at the University of Basel. Paul Euler had high hopes that, following in his footsteps, his son would become a Protestant minister. However, it was not long before Johann, who was Leonhard's mentor, noticed the young boy's mathematical ability while he was a student (at the age of fourteen) at the University of Basel. In Euler's own words:

> I soon found an opportunity to be introduced to a famous professor Johann Bernoulli. ... True, he was very busy and so refused flatly to give me private lessons; but he gave me much more valuable advice to start reading more difficult mathematical books on my own and to study them as diligently as I could; if I came across some obstacle or difficulty, I was given permission to visit him freely every Sunday afternoon and he kindly explained to me everything I could not understand…[18]

Given his close relationship with the Bernoullis, it is not surprising that Euler became interested in the calculus of variations. As early as 1728, Leonhard Euler had already written "On finding the equation of geodesic curves." By the 1730s, he was concerning himself with isoperimetric problems.

In 1744, Euler published his landmark book *Methodus inveniendi lineas curvas maximi minimive proprietate gaudentes, sive solutio problematis isoperimetrici latissimo sensu* *accepti* (*A method for discovering curved lines that enjoy a maximum or minimum property, or the solution of the isoperimetric problem taken in the widest sense*). Some mathematicians date this as the birth of the *theory* of the calculus of variations [14].

Euler took the methods used to solve specific problems and systematized them into a powerful apparatus. With this method, he was then able to study a very general class of problems. His opus considered a variety of geodesic problems, various modified and more general brachistochrone problems (such as considering the effects of a resistance to the falling body), problems involving isoperimetric constraints, and even questions of invariance. Although few mathematicians before Euler would give a second thought to such things, he examined whether his fundamental conditions would remain unchanged under general coordinate transformations. (These questions were not completely resolved until the twentieth century.)

Also in this publication, it was shown for the first time that in order for $y(x)$, satisfying

$$I = \int_{x_A}^{x_B} f(x, y, y')dx, \ y(x_A) = y_A, \ y(x_B) = y_B, \ x_A < x_B,$$

to yield a minimum of *I*, then a necessary condition is the so-called Euler-Lagrange equation (which first appeared in Euler's work eight years previously)

$$\frac{\partial f}{\partial y} - \frac{d}{dx}\left(\frac{\partial f}{\partial y'}\right) = 0.$$

Another important element of Euler's exposition was his statement and discussion of a very important principle in mechanics. However, it has also been attributed to another, lesser, mathematician.

In two papers read to l'Académie des Sciences in 1744, and to the Prussian Academy in 1746, the French mathematician Maupertuis (1698-1759) proclaimed to the world *le principe de la moindre quantité d'action*, or the principle of least action. In an almost Pythagorean spirit, Maupertuis said that "*la Nature, dans la production de ses effets, agit toujours par les moyens les plus simples*." (In her actions, Nature always works by the simplest methods.) He believed that physically, things unfold in Nature in such a way that a certain quality, which he called the "action," is always minimized. While Maupertuis's intuition was good, he certainly lacked the logical motivation and clarity of Euler. His definition of the action was vague. His rationale in developing this principle was somewhat mystical. He sought to develop not only a mathematical foundation for mechanics but a theological one as well. He went so far as to say, in his *Essai de cosmologie* (1759), that the perfection of God would be incompatible with anything other than utter simplicity and the minimum expenditure of *action*!

> *Notre principe, plus conforme aux idées que nous devons avoir des choses, laisse le Monde dans le besoin continuel de la puissance de Créateur, et est une suite nécessaire de l'emploi le plus sage de cette puissance… Ces loix si belle et si simples sont peut-être les seules que le Créateur et l'Ordonnateur des choses a établies dans la matière pour y opérer tous les phénomènes de ce Monde visible.*

> (Our principle, which conforms better to the ideas that we should have about things, leaves the world inconstant need of the strength of the Creator and follows necessarily from the most wise use of this strength… These simple and beautiful laws are perhaps the only ones that the Creator and Organizer of all things has put in place to carry out the workings of the visible world.) [20]

Returning to Euler, and his magnificent work of 1744, we see strikingly similar ideas but without the theological overtones. Near the beginning of the section on the principle of least action, Euler writes:

> Since all the effects of Nature follow a certain law of maxima or minima, there is no doubt that, on the curved paths, which the bodies describe under the action of certain forces, some maximum or minimum property ought to obtain. What this property is, nevertheless, does not appear easy to define *a priori* by proceeding from the principles of metaphysics; but since it may be possible to determine these same curved paths by means of a direct method, that very thing which is a maximum or minimum along these curves can be obtained with due attention being exhibited. But above all the effect arising from the disturbing forces ought especially to be regarded; since this [effect] consists of the motion produced in the body, it is consonant with the truth that this same motion or rather the aggregate of all motions, which are present in the body ought to be a minimum. Although this conclusion does not seem sufficiently confirmed, nevertheless if I show that it agrees with a truth known *a priori* so much weight will result that all doubts which could originate on this subject will completely vanish. Even better when its truth will have been shown, it will be very easy to undertake studies in the profound laws of Nature and their final causes, and to corroborate this with the firmest arguments [11].

As often happens in mathematics even today, there was a bitter dispute as to the priority of the discovery of the principle of least action. In 1757, the mathematician König produced a letter supposedly written by Leibniz in 1707 that contained a formulation of the principle of least action. At the time, Maupertuis, who was a headstrong and virulent man, was the president of the Prussian Academy and had a sharp reaction to this claim. He accused his fellow-member of plagiarism and was convinced that the letter was a forgery. Ironically, Euler sided with his French colleague in this affair, even though it is possible (and perhaps most likely) that it was Euler himself who was the first to put his finger on the principle.

An additional topic of interest stemming from Euler's opus of 1744 is that of minimal surfaces. One of the most fascinating areas of geometry, minimal surfaces are obtained from the calculus of variations as portions of surfaces of least area among all surfaces bounded by a given space curve. Euler discovered the first non-trivial such surface, the catenoid, which is generated by rotating a catenary (i.e. a cosh curve or the curve of a hanging chain); for example, $r = A \cosh x$, where r is the distance in 3-dimensional space from the x-axis [5]. We will have more to say about minimal surfaces later.

While it is true that a short time later, Euler's technique was superseded by that of Lagrange (as we shall soon see), at the time it was completely new mathematics. His systematic methods, in an elegant form, were remarkable for their clarity and insight. As the twentieth century mathematician Carathéodory, who edited Euler's works, wrote in the introduction,

> [Euler's book] is one of the most beautiful mathematical works ever written. We cannot emphasize enough the extent to which that Lehrbuch over and over again served later generations as a prototype in the endeavour of presenting special mathematical material in its [logical, intrinsic] connection [14].

**Lagrange**

In 1755, a 19-year-old from Turin sent a letter to Euler that contained details of a new and beautiful idea. Euler's correspondent, Ludovico de la Grange Tournier, was no ordinary teenager. Less than two months after he wrote that fateful letter to Euler, the man we now know as Joseph-Louis Lagrange (1736-1813) was appointed professor of mathematics at the Royal Artillery School in Turin. His rare gifts, his humility, and his devotion to mathematics made him one of the giants of eighteenth century mathematics. He contributed much groundbreaking work in fields as diverse analysis, number theory, algebra, and celestial mechanics. However, it was with the calculus of variations that his early reputation was made.

In his first letter to the legendary Swiss mathematician, Lagrange showed Euler how to eliminate the tedious geometrical methods from his process. Essentially, Lagrange had developed the idea of comparison functions (like the $\eta(x)$ function used in the mathematical background section above), which lead almost directly to the Euler-Lagrange equation. After considering Lagrange's method, Euler became an instant convert, dropped his old geometrical methods, and christened the entire field by the name we now use, the calculus of variations, in honour of Lagrange's variational method [11].

With the recipe reduced to a much simpler analytic method, even more general results could be obtained. The following year, in 1756, Euler read two papers to the Berlin Academy in which he made liberal use of Lagrange's method. In his first paper, he was quick to give the young man from Turin his due:

> Even though the author of this [Euler] had meditated a long time and had revealed to friends his desire yet the glory of first discovery was reserved to the very penetrating geometer of Turin, Lagrange, who having used analysis alone, has clearly attained the very same solution which the author had deduced by geometrical considerations [11].

The two great mathematicians corresponded frequently over the next few years, with Lagrange working hard to extend the theory. Toward the end of 1760, he was able to publish a number of his results in *Miscellanea Taurinensia*, a scientific journal in Turin, under the title *Essai d'une nouvelle méthode pour déterminer les maxima et les minima des formules intégrals indefinites* (*Essay on a new method for determining maxima and* *minima for formulas of indefinite integrals*).  Solutions to more general problems we investigated for the first time, such as variable end-point brachistochrone problems, finding the surface of least area among all those bounded by a given curve (a problem that we associate today with Plateau), and finding the polygon whose area is greatest among all those that have a fixed number of sides.  An apt résumé of the advances of the new theory comes from the pen of Lagrange himself:

> Euler is the first who has given the general formula for finding the curve along which a given integral expression has its greatest value…but the formulas of this author are less general than ours: 1. because he only permits the single variable *y* in the integrand to vary; 2. because he assumes that the first and last points of the curve are fixed…By the methods which have been explained one can seek the maxima and minima for curved surfaces in a most general manner that has not been done till now [11].

It was also in this early work of Lagrange that his famous rule of multipliers was first discussed.  However, the generality and power of the method was not clear to him at that time and it was not until his path-breaking tour de force *Méchanique analytique* (1788), that he clearly expressed the rule in its modern form.

When trying to extremize a function, often difficulties arise when the function is subject to certain outside conditions or constraints.  In principle, we could use each constraint to eliminate one variable at a time, thereby reducing the problem progressively to a simpler and simpler one.  However, this can be both tedious and time consuming.  Lagrange's method of multipliers is a powerful tool that allows for solutions to the problem without having to solve the conditions or constraints explicitly.  Let us now show the solution of such a problem, arising from a simple quantum mechanical system.

> Consider the problem of a particle of mass *m* in a box, which we can consider as a parallelepiped with sides *a*, *b*, and *c*.  The so-called ground state energy of the particle is given by
>
> $$E = \frac{h^2}{8m}\left(\frac{1}{a^2} + \frac{1}{b^2} + \frac{1}{c^2}\right),$$
>
> where *h* is Planck's constant.  Now suppose we wish to find the shape of the box that will minimize the energy *E*, subject to the constraint that the volume of the box is constant, i.e.
>
> $$V(a,b,c) = abc = k.$$

Essentially, we need to minimize the function $E(a,b,c)$ subject to the constraint $\varphi(a,b,c) = abc - k = 0$.  For the variable *a*, this implies that

$$\frac{\partial E}{\partial a} + \lambda\frac{\partial \varphi}{\partial a} = -\frac{h^2}{4ma^3} + \lambda bc = 0,$$

where $\lambda$ is an arbitrary constant (called the Lagrange multiplier). Of course we have similar equations for the other variables:

$$-\frac{h^2}{4mb^3} + \lambda ac = 0, \qquad -\frac{h^2}{4mc^3} + \lambda ab = 0.$$

After multiplying the first equation by *a*, the second by *b*, and the third by *c*, we obtain

$$\lambda abc = \frac{h^2}{4ma^2} = \frac{h^2}{4mb^2} = \frac{h^2}{4mc^2}.$$

Hence, our solution is $a = b = c$, which is a cube. Notice how we did not even need to determine the multiplier $\lambda$ explicitly [2].

The *Méchanique analytique* was an ambitious undertaking, as it summarized all the work done in the field of classical mechanics since Newton. In fact, as books on mechanics go, it is mentioned in the same breath as Newton's *Philosophiae naturalis principia mathematica*. Whereas Newton considered most problems from the geometrical point of view, Lagrange did everything with differential equations. In the preface, he even states that

> …one will not find figures in this work. The methods that I expound require neither constructions, nor geometrical or mechanical arguments, but only algebraic operations, subject to a regular and uniform course [11].

Classical mechanics had really come of age with Lagrange. Building on the great insights of Euler, Lagrange was able to rescue mechanical problems from the tedium of geometrical methods. His approach is still meaningful today and it forms one of the cornerstones of the mathematical framework of modern theoretical physics. As it turns out, there was still much work to be done in the calculus of variations. There were unforeseen problems with the approach of Euler and Lagrange. However, let us pay our debt to Lagrange by remembering the words of Carl Gustav Jacob Jacobi (1804-1851), who was one the main contributors to the theory of variational problems in the nineteenth century:

> By generalizing Euler's method he arrived at his remarkable formulas which in one line contain the solution of all problems of analytical mechanics.
>
> [In his Memoir of 1760-61] he created the whole calculus of variations with one stroke. This is one of the most beautiful articles that has ever been written. The ideas follow one another like lightning with the greatest rapidity [14].

**Legendre**

In 1786, Adrien-Marie Legendre (1752-1833) presented a memoir to the Paris Academy entitled *Sur la manière de distinguer les maxima des minima dans le calcul des variations* (*On the method of distinguishing maxima from minima in the calculus of variations*). Legendre was a well-known mathematician from Paris who developed many analytical tools for problems in mathematical physics and served as editor for Lagrange's *Méchanique analytique*.

Legendre considered the problem of determining whether an extremal is a minimizing or a maximizing arc. Let us recall that in extrema problems of one variable calculus, we consider not only points where the first derivative vanishes, but we also study the second derivative at these points. Similarly, Legendre examined the "second variation" of the functional, motivated by the theorem of Taylor:

$$\delta^2 I = \frac{\varepsilon^2}{2} \frac{\partial^2 I[\tilde{y}]}{\partial \varepsilon^2}\bigg|_{\varepsilon=0} = \int_{x_A}^{x_B} \left( \frac{\varepsilon^2}{2} f_{yy} \eta^2 + 2 f_{yy'} \eta \eta' + f_{y'y'} \eta'^2 \right) dx.$$

Legendre was able to show the condition $f_{y'y'} \geq 0$ along a minimizing curve and $f_{y'y'} \leq 0$ along a maximizing curve, which is surprisingly similar to what we obtain in elementary calculus in the second derivative test! In spite of the fact that he was on the right track, Legendre's attempt to show that this condition is both necessary and sufficient was not quite correct [11], [14]. The idea did not catch on and by the time Lagrange levelled several objections to the second variation approach in his *Théorie des fonctions analytiques* (1797), it appeared that the death knell has sounded for Legendre's innovative idea.

**Jacobi**

It was not until fifty years passed since Legendre's initial discovery of the second variation condition that another mathematician took up the task of developing the theory even further. In 1836, in a paper remarkable for its brevity and obscurity, Jacobi demonstrated rigourously what we now call the Jacobi condition, namely that:

> For a local minimum, it suffices to have both of the following satisfied:

1) $f_{y'y'} > 0$, and
2) $x_B$ closer to $x_A$ than to the "conjugate point" of $x_A$, which is the first value $x > x_A$ where a nonzero solution of

$$\frac{d}{dx}\left( f_{y'y'} \frac{dw}{dx} \right) - \left( f_{yy} - \frac{d}{dx} f_{yy'} \right) w = 0, \quad w(x_{x_A}) = 0 \quad (x \geq x_A)$$

vanishes. [Here $w(x) = \left.\frac{\partial y}{\partial \alpha}\right|_{\alpha=0}$ and $\alpha = 0$ corresponds to $\tilde{y}$ in the family of extremals $y = y(x, \alpha)$.]

Another way to say this is that when the two conditions above are satisfied, then there exists a minimizing $\tilde{y}$ among $y \in C^1[x_A, x_B]$ satisfying the boundary conditions $y(x_A) = y_A$ and $y(x_B) = y_B$, and satisfying:

(a) $|y - \tilde{y}| < \rho,$  (b) $|y' - \tilde{y}'| < \rho$ for small positive $\rho$.

This paper, entitled *Zur Theorie der Variations-Rechnung und der Differential-Gleichungen* (*On the calculus of variations and the theory of differential equations*), was so terse that rigourous proofs were not given but instead were hinted at [11], [14]. Perhaps, as one mathematical historian has suggested, Jacobi was in a rush to publish his results first to ensure intellectual priority [11]. It is difficult to agree with such a theory since progress in this field had stagnated for half a century! In any case, it was an opportunity for numerous mathematicians to provide further elucidation and commentary in the years that followed.

**Hamilton-Jacobi Theory**

While not directly connected with the development of the theory of the calculus of variations, it is timely to draw attention to another aspect of Jacobi's work. In the mid 1830s, a Scottish mathematician named William Rowan Hamilton (1788-1856) developed the foundations of what we now call Hamiltonian mechanics. Closely related to the methods developed by Lagrange, Hamilton showed that under certain conditions, problems in mechanics involving many variables and constraints can be reduced to an examination of the partial derivatives of a single function, which we now appropriately call the Hamiltonian. In the original papers of 1834 and 1835, some rigour was lacking and Jacobi was quick to step in. Hamilton did not show under which conditions he could be certain that his equations possessed solutions. In 1838, Jacobi was able to rectify this, in addition to showing that one of the equations Hamilton studied was redundant. Due to the tidying up and simplification performed by Jacobi, many modern books on classical mechanics refer to this approach as Hamilton-Jacobi theory [11].

**Nineteenth Century Applications to Other Fields: Edgeworth and Poe**

By the nineteenth century, mathematical methods had advanced further than many had dreamed possible. Previously unsolved problems in physics, astronomy, engineering, and technology were being overcome at last. New theories were being developed at a speed never seen before, with a startling predictive nature that few imagined possible. One only needs to consider Newtonian mechanics, the developments in understanding thermodynamic systems, or especially, the elegant systematization of the theories of electricity and magnetism laid out in Maxwell's equations. How natural, then that people tried to apply the same powerful techniques to other disciplines. In some cases, a measure of success was attained. In other cases, the results seem laughable.

In 1881, a book appeared with the title *Mathematical Psychics: An Essay on the Application of Mathematics to the Moral Sciences* [19]. The author was Francis Edgeworth (1845-1926), an English economist. A primary goal of the text was to construct a model of human science in which ethics can be viewed as a science. Today, the book is remembered chiefly for the merit of its ideas for economic theory. For us, the most interesting part of the book is the section on utilitarian calculus. Inspired by the utilitarian Jeremy Bentham (1748-1832), Edgeworth used the mathematical techniques of the calculus of variations in an effort to extremize the happiness function, or a function that was designed to measure the achievement of the ultimate good in society.

Defining fundamental units of pleasure within the context of human interpersonal contracts, Edgeworth was able to obtain an equation involving the sum over all individuals' utility. Despite variations from point to point, Edgeworth hypothesized that there would exist a locus at which the sum of the utilities of the individuals is a maximum. Edgeworth called this the *utilitarian point*. Edgeworth was quick to realize that the Benthamite slogan, "the greatest happiness of the greatest number" needed restating in a more precise form. After some mathematical labour, he was able to show that "the ultimate good was to be conceived as the maximum value of the triple integral over the variables 'pleasure,' individuals, and time."

In retrospect, it is hardly surprising that this treatise has no impact on the development of moral and ethical philosophy.

Caught up in the spirit of things, and inspired by the writings of the greatest mathematicians on the calculus of variations, Edgar Allan Poe (1809-1849) published a story in 1841 called *Descent into the Maelstrom* [12]. In the story, the protagonist is able to survive a violent storm by noting certain critical properties of solids moving in a resisting medium:

> ...what I observed was, in fact, the natural consequence of the forms of floating fragments...a cylinder, swimming in a vortex, offered more resistance to its suction, and was drawn in with greater difficulty than any equally bulky body, of any form whatever.

Poe was inspired, no doubt, by Newton's Principia. Fortunately for Poe, good science is not needed in order to tell a good story. In the story, it is claimed that the sphere offered the minimum resistance, although Newton showed long ago that this is not the case. In addition, Newton's results were only good for bodies moving through a motionless fluid, not a violent sea. In any case, it is still a good example of how science can motivate the creative arts.

**Riemann, Dirichlet, and Weierstrass**

It is surprising to discover that the development of the theory of the calculus of variations not only impacted physical problems and the theory of partial differential equations, but also the fields of classical analysis and functional analysis. In the mid-1800s, many mathematicians, such as Bernhard Riemann (1826-1866) and Gustave Lejeune Dirichlet (1805-1859) searched for general solutions to boundary value and initial value problems of partial differential equations arising in physical problems. Problems of this type are of great importance in physics, as they are basic to the understanding of gravitation, electrostatics, heat conduction, and fluid flow. One of the problems that attracted many of the top mathematicians of the day was an existence proof of a solution $u$, in a general domain $\Omega$, satisfying:

$$\nabla^2 u = 0 \text{ in } \Omega; \quad u\big|_{\partial\Omega} = f, \quad u \in C^2(\Omega) \cap C^0(\overline{\Omega}), \quad \Omega \subset \mathrm{R}^2 \text{ or } \mathrm{R}^3,$$

where $\nabla^2 u = u_{xx} + u_{yy} + u_{zz}$. This is known as a Dirichlet problem. Riemann used principles from the calculus of variations to develop a proof of this, which was a problem he had first seen in lectures by Dirichlet. He named it Dirichlet's principle and stated it as follows

> There exists a function $u$ that satisfies the condition above and that minimizes the functional
>
> $$D[u] = \int_\Omega |\nabla u|^2 \, dV, \qquad \Omega \subset \mathrm{R}^2 \text{ or } \mathrm{R}^3,$$
>
> among all functions $u \in C^2(\Omega) \cap C^0(\overline{\Omega})$ which take on given values f on the boundary $\partial\Omega$ of $\Omega$.

Dirichlet's principle had been used earlier by Gauss (1839) and Lord Kelvin (1847) before Riemann used the principle in 1851 in order to obtain fundamental results in potential theory using complex analytic functions [15], [16]. However, something was not quite right with the theory. As one mathematician noted:

> It was a strange situation. Dirichlet's principle had helped to produce exciting basic results but doubts about its validity began to appear, first in private remarks of Weierstrass - which did not impress Riemann, who placed no decisive value on the derivation of his existence theorems by Dirichlet's principle - and then, after both Dirichlet and Riemann had died, in Weierstrass's public address to the Berlin Academy...

As it turns out, there was a fundamental conceptual error involved in the faulty method of proof employed by Riemann. He failed to distinguish the differences between a greatest lower bound and a minimum for the Dirichlet problem. Karl Weierstrass (1815-1897) was the first to point out that in some cases, a minimizing function can come arbitrarily close to the lower bound without ever reaching it.

The breakdown of Dirichlet's principle (which had been the basis for many new results) turned out to be very beneficial for the theory of analysis. In an effort to patch up the theory, three new methods of existence proofs were developed, by Hermann Schwarz (1843-1921), Henri Poincaré (1854-1912), and Carl Neumann (1832-1925) [15].

Beginning in the 1870s, Weierstrass gave the theory of the calculus of variations a complete overhaul. It took quite some time for these results to become widely known to the rest of the mathematical community, principally through the dissertations of his graduate students. Known for his rigourous approach to mathematics, Weierstrass was the first to stress the importance of the domain of the functional that one is trying to minimize. He also examined the family of admissible functions satisfying all of the constraints. His most notable accomplishment was the fact that he gave the first ever completely correct sufficiency theorem for a minimum. Two new concepts, the field of extremals and the E-function, were developed in order to tackle the problem of sufficiency and a new type of minimum (a so-called strong minimum) was defined [15], [16].

**Philosophical Interlude**

To the applied mathematician or physicist, all of this work to define conditions of sufficiency for the existence of an extremum might sound like splitting hairs. As Göthe wrote in *Maxims and Reflections*,

> Mathematicians are like a certain type of Frenchman: when you talk to them they translate it into their own language, and then it soon turns into something completely different.

For problems in mechanics, for example, the Euler-Lagrange equation works perfectly well ninety-nine times out of a hundred - and when it doesn't, then it should be physically obvious. This point of view was expressed by Gelfand and Fomin:

> ...the existence of an extremum is often clear from the physical or geometric meaning of the problem, e.g., in the brachistochrone problem... If in such a case there exists only one extremal satisfying the boundary conditions of the problem, this extremal must perforce be the curve for which the extremum is achieved [17].

The rigourous mathematician would surely answer that in mathematics, conclusions should be logically deducible from initial hypotheses. And when it comes to a physical model, the mathematician would no doubt remind us that we should be mindful of the assumptions and idealizations we make for the sake of simplicity, and the consequences these assumptions entail.

In reality, what is truly surprising is not that mathematicians fought over the smallest details of the calculus of variations for more than one hundred years, but that it took so long for anyone to realize the elementary mistakes that Euler made when he first examined these problems. A twentieth century mathematician, L.C. Young, remarked at length on this oversight in his excellent book, *Lectures on the Calculus of Variations and Optimal Control Theory* [21]. It is rewarding to see how he puts things into perspective:

> In the Middle Ages, an important part was played by the jester: a little joke that seemed so harmless could, as its real meaning began to sink in, topple kingdoms. It is just such little jokes that play havoc today with a mathematical theory: we call them paradoxes.
>
> Perron's paradox runs as follows: "Let $N$ be the largest positive integer. Then for $N \neq 1$ we have $N^2 > N$ contrary to the definition of $N$ as largest. Therefore $N = 1$."
>
> The implications of this paradox are devastating. In seeking the solution to a problem, we can no longer assume that this solution exists. Yet this assumption has been made from time immemorial, right back in the beginnings of elementary algebra, where problems are solved starting off with the phrase: "Let *x* be the desired quantity."
>
> In the calculus of variations, the Euler equation and the transversality conditions are among the so-called necessary conditions. They are derived by exactly the same pattern of argument as in Perron's paradox; they assume the existence of a solution. This basic assumption is made explicitly, and it is then used to calculate the solutions whose existence was postulated. In the class of problems in which the basic assumption is valid, there is nothing wrong with doing this. But what precisely *is* this class of problems? How do we know that a particular problem belongs to this class? The so-called necessary conditions do not answer this. Therefore a "solution" derived by necessary conditions only is simply no valid solution at all.
>
> It is strange that so elementary a point of logic should have passed unnoticed for so long! The first to criticize the Euler-Lagrange method was Weierstrass, almost a century later. Even Riemann made the same unjustified assumption in his famous Dirichlet principle...
>
> The main trouble is that, as Perron's paradox shows, the fact that a "solution" has actually been calculated in no way disposes of the logical objection to the original assumption.
>
> A reader may here interpose that, in practice, surely this is not serious and would lead no half competent person to false results; was not Euler at times logically incorrect by today's standards, but nonetheless correct in his actual conclusions? Do not the necessary corrections amount to no more than a sprinkling of definitions, which his insight perhaps took into account, without explicit formulation?
>
> Actually, this legend of infallibility applies neither to the greatest mathematicians nor to competent or half competent persons, and the young candidate with an error in his thesis does not disgrace his calling... Newton formulated a variational problem of a solid of revolution of least resistance, in which the law of resistance assumed is physically absurd and ensures that the problem has no solution – the more jagged the profile, the less the assumed resistance – and this is close to Perron's paradox. If this had been even approximately correct, after removing absurdities, there would be no need today for costly wind tunnel experiments. Lagrange made many mistakes. Cauchy made one tragic error of judgment in rejecting Galois's work. The list is long. Greatness is not measured negatively, by absence of error, but by methods and concepts which guide further generations [21].

**Twentieth Century Developments**

With the calculus of variations on a relatively firm foothold, aided by the rigourous work of the school of Weierstrass, things were set for the theory to develop even further. In his famous turn- of-the-century address to the International Congress of Mathematicians in Paris, David Hilbert (1862-1943) made mention of the calculus of variations on several occasions when discussing other problems. In addition, his twenty-third problem was a call for the further elucidation of the theory:

> So far, I have generally mentioned problems as definite and special as possible, in the opinion that it is just such definite and special problems that attract us the most and from which the most lasting influence is often exerted upon science. Nevertheless, I should like to close with a general problem, namely with the indication of a branch of mathematics repeatedly mentioned in this lecture— which, in spite of the considerable advancement lately given it by Weierstrass, does not receive the general appreciation which, in my opinion, is its due—I mean the calculus of variations.

In the next few years, Hilbert and his associates continued where Weierstrass left off, developing many new results and setting the stage for the next leap forward.

**Morse Theory**

Marston Morse (1892-1977) turned his eye to the global picture and developed the calculus of variation in the large, with applications to equilibrium problems in mathematical physics. We now call the field Morse theory. In a paper published in 1925 entitled *Relations between the critical points of a real function of n independent variables*, Morse proved some important new results that had a big effect on global analysis, which is the study of ordinary and partial differential equations from a topological point of view. Much of his work depended on the results obtained by Hilbert and company [15].

**Optimal Control Theory**

Another new field developed in the twentieth century from the roots of the calculus of variations is optimal control theory. A generalization of the calculus of variations, this theory is able to tackle problems of even greater generality and abstraction. New mathematical tools were developed by chiefly Pontryagin, Rockafellar, and Clarke that, among other things, enabled nonlinear and nonsmooth functionals to be optimized. While this may sound like a mathematical abstraction, in reality there are many physical problems that can only be solved in such a manner. Two examples which come from the engineering world are the problem of landing a spacecraft as softly as possible with the minimum expenditure of fuel and the construction an ideal column [9].

**Minimal Surfaces**

The minimal surfaces discovered by Euler have also played a substantial role in twentieth century mathematics, during which time two Fields Medals were awarded for work related to the subject. In 1936, Jesse Douglas won a Fields Medal for his solution to Plateau's problem and in 1974, Enrico Bombieri shared a Fields Medal for his work on higher dimensional minimal surfaces. It is becoming apparent that minimal surfaces are found throughout nature. Examples are soap films, grain boundaries in metals, microscopic sea animals (called radiolarians), and the spreading and sorting of embryonic tissues and cells. In addition, minimal surfaces have proved popular in design, through the work of the German architect Frei Otto, as well as in art, exemplified in the works of J.C.C. Nitsche [1].

**Physics**

We have already seen the rich interplay between the mathematical methods used in the calculus of variations and developments in understanding the natural laws of our universe. Recall the least time principles of Fermat, Maupertuis, Euler, Lagrange, and Hamilton and their effects on the history of optics and mechanics. The success of these variational methods in solving physical problems is not surprising [9]. As Yourgrau and Mandelstam point out:

> Arguments involving the principle of least action have excited the imagination of physicists for diverse reasons. Above all, its comprehensiveness has appealed, in various degrees, to prominent investigators, since a wide range of phenomena can be encompassed by laws differing in detail yet structurally identical. It seems inevitable that some theorists would elevate these laws to the status of a single, universal canon, and regard the individual theorems as mere instances thereof. It further constitutes an essential characteristic of action principles that they describe the change of a system in such a manner as to include its states during a definite time interval, instead of determining the changes which take place in an infinitesimal element of time, as do most differential equations of physics. On this account, variational conditions are often termed "integral" principles as opposed to the usual "differential" principles. By enforcing seemingly logical conclusions upon arguments of this type, it has been claimed that the motion of the system during the whole of the time interval is predetermined at the beginning, and thus teleological reflections have intruded into the subject matter. To illustrate this attitude: if a particle moves from one point to another, it must, so to speak, 'consider' all the possible paths between the two points and 'select' that which satisfies the action condition [20].

In 1948, motivated by a suggestion by P.A.M. Dirac, the American physicist Richard Feynman (1918-1988) developed a completely new approach to quantum mechanics, based on variational methods. Although not mathematically well-defined, the Feynman path integral was what he called a "summation over histories" of the path of a particle. Despite the fact that the original paper was rejected by one journal for being nothing new, Feynman's original approach was ideally suited to extending quantum theory to a more general framework, incorporating relativistic effects [10].

It did not take long for the mathematicians to come along and tidy up everything. Mark Kac showed that Feyman's integral can be thought of as a special case of the Wiener integral, developed by Norbert Wiener in the 1920s. With a rigorous mathematical underpinning, physicists were then able to apply the new variational techniques to a host of all quantum and statistical phenomena. Today, these methods are employed in the monumental task of developing the so-called Grand Unified Theory.

As the field evolved from our search to understand the inner workings of Nature, perhaps it is fitting to end this survey of the history of the calculus of variations with a quote from someone still actively involved in this search. When asked about the role of the calculus of variations in modern physics, Maxim Pospelov, a theoretical physicist specializing in supersymmetry, had this to say:

> The most notable change that the 20th century brought to physics is the transition from a deterministic classical mechanics where the variation of action leads to the equations of motion and single trajectory when the boundary conditions are fixed to quantum mechanics that allows multiple trajectories and determines the probability for a certain trajectory. The functional integral approach to quantum mechanics and quantum field theory is the modern language that everybody uses. All, absolutely all, physical processes in quantum field theory can be studied as a variation of the vacuum-vacuum transition amplitude in the presence of external sources over these sources.
>
> Variational methods are often used in particular calculations when, for example, one needs to find a complicated wave function when the exact solution to the Schrödinger equation is not possible. I know that the variational approach to the helium atom yields a very precise determination of its energy levels and ionization threshold [7].

**Bibliography**

[1]   Almgren F.J. (1982)  Minimal surface forms.  *Math. Intelligencer* Vol. 4 No.4, pp. 164-172.

[2]   Arfken G. and Weber H. (2001)  *Mathematical Methods for Physicists*. San Diego: Academic Press.

[3]   Ball J.M. (1998)  The calculus of variations and materials science.  *Quart. Appl. Math.* Vol. 56, No. 4, pp. 719-740.

[4]   Buttazzo G. and Kawohl B. (1993)  On Newton's problem of minimal resistance.  *Math. Intelligencer* Vol. 15 No.4, pp. 7-12.

[5]   Byron F. and Fuller R. (1969) *Mathematics of Classical and Quantum Physics*.  Reading: Addison-Wesley.

[6]   Cuomo S. (2000)  *Pappus of Alexandria and the Mathematics of Late Antiquity*.  Cambridge: Cambridge University Press.

[7]   Ferguson J. (2003) Private e-mail correspondence with M. Pospelov.

[8]   Ferguson J. (2003) Private discussions with J. Ye.

[9]   Ferguson J. (2003) Private discussions with W. Israel.

[10] Feynman, R. (1948) Space-time approach to non-relativistic quantum mechanics.  *Rev. Mod. Phys.* Vol. 20, No. 2, pp.367-387.

[11] Goldstine H. (1980)  *A History of the Calculus of Variations from the 17th through the 19th Century*. New York: Springer-Verlag.

[12] Gould S. (1985)  Newton, Euler, and Poe in the calculus of variations.  *Differential geometry, calculus of variations, and their applications.*  Gould S. (Ed.) New York: Dekker.

[13] Kirmser P. and Hu K-K. (1993)  The shape of an ideal column reconsidered.  *Math. Intelligencer* Vol. 15 No.3, pp. 62-68.

[14] Kreyszig E.  On the calculus of variations and its major influences on the mathematics of the first half of our century. I.  *Amer. Math. Monthly* 101 (1994) no.7, pp. 674-678.
[15] _________.  On the calculus of variations and its major influences on the mathematics of the first half of our century. II.  *Amer. Math. Monthly* 101 (1994) no.9, pp. 902-908.
[16] _________.  Interaction between general topology and functional analysis.  *Handbook of the History of General Topology, Vol.1*, Aull C.E. and Lowen R. (Eds.), pp. 357-389, Kluwer Acad. Publ., Dordrecht, 1997.

[17] McShane E.J. (1989)  The calculus of variations from the beginning through optimal control theory.  *SIAM J. Cont. Optim.*  Vol. 27, No. 5, pp. 916-989

[18] O'Connor J.J. and Robertson E.F. *MacTutor History of Mathematics Archive*. http://www-gap.dcs.st-and.ac.uk/~history/. 29 Nov. 2003.

[19] Wall B. (1978/79) F. Y. Edgeworth's mathematical ethics. Greatest happiness with the calculus of variations. *Math. Intelligencer* Vol. 1, No.3, pp. 177-181.

[20] Yourgrau W. and Mandelstam S. (1968) *Variational Principles in Dynamics and Quantum Theory*. London: Pitman & Sons.

[21] Young L.C. (1969) *Lectures on the Calculus of Variations and Optimal Control Theory*. Philadelphia: W.B. Saunders Company.